\documentclass[11pt,a4paper]{article}
\usepackage{amssymb,  amsthm, color, amsfonts, amsmath}
\usepackage{graphics}
\usepackage{latexsym}
\title{\vspace{-2cm}
Low Mach number limit for the Quantum-Hydrodynamics system
}
\author{\textbf{Donatella Donatelli}\\
        {\small Department of Information Engineering, Computer Science and Mathematics}\\
       {\small University of L'Aquila}\\
       {\small 67100 L'Aquila, Italy}\\
        $\scriptstyle\mathtt{donatella.donatelli@univaq.it}$\\
 \and        
  \textbf{ Pierangelo Marcati}\\
        {\small Department of Information Engineering, Computer Science and Mathematics}\\
       {\small University of L'Aquila}\\
       {\small and GSSI - Gran Sasso Science Institute}\\
       {\small 67100 L'Aquila, Italy}\\
        $\scriptstyle\mathtt{pierangelo.marcati@univaq.it,\ pierangelo.marcati@gssi.infn.it}$\\}

         \date{}

\newcommand{\e}{\varepsilon}		       
\newcommand{\R}{\mathbb{R}}

\newcommand{\T}{\mathbb{T}}

\newcommand{\Je}{J^{\varepsilon}}
\newcommand{\Le}{\Lambda^{\varepsilon}}

\newcommand{\re}{\rho^{\varepsilon}}

\newcommand{\se}{\delta^{\varepsilon}}

\newcommand{\dive}{\mathop{\mathrm {div}}}

\newtheorem{theorem}{Theorem}[section]

\newtheorem{proposition}[theorem]{Proposition}

\theoremstyle{definition}
\newtheorem{definition}[theorem]{Definition}
\newtheorem{remark}[theorem]{Remark}

\begin{document}
\maketitle
\begin{abstract}
In this paper we deal with the low Mach number limit for the system of  quantum-hydrodynamics, far from the vortex nucleation regime. More precisely, in the framework of a periodic domain and ill-prepared initial data we prove strong convergence of the solutions towards regular solutions of the incompressible Euler system. In particular we  will perform a detailed analysis of the time oscillations and of the relative entropy functional related to the system.  

\medbreak 
\medbreak
\end{abstract}
\newpage
\tableofcontents
\newpage
\section{Introduction}
In this paper we deal with the following system of  quantum hydrodynamics 
\begin{equation}
\partial_{s}\rho+\dive J=0,
\label{1}
\end{equation}
\begin{equation}
\partial_{s} J+\dive\left(\frac{J\otimes J}{\rho}\right)+\nabla p(\rho)=\dive(\rho\nabla^{2}\log \rho),
\label{2}
\end{equation}
where $\rho$  and $J$ represent the charge and current density respectively and $p(\rho)$ is the hydrodynamic pressure which is a function depending only of $\rho$, satisfying  the following conditions
\begin{equation} 
p(\rho) \in  C[0,\infty) \cap C^1(0,\infty) \label{2.1}
\end{equation}
and
\begin{equation}
\begin{cases}
p(0) = 0, \,\, p'(\rho) \ge a_1 \rho^{\gamma-1}-b &  \mbox{for all  $\rho >0$}, \\
p(\rho) \le a_2 \rho^{\gamma} +b, & \mbox{for all $\rho \geq 0$, $\gamma \geq \frac{n}{2}$}. 
\end{cases}
\label{3} 
\end{equation}
The computations we are going to perform later on  can be easily adapted to the general pressure law \eqref{3} however for simplicity in this paper we will take $p(\rho)$ of the form
\begin{equation}
p(\rho)=\rho^{\gamma}/\gamma.
\label{4}
\end{equation}
By using the Madelung formalism the quantum hydrodynamics with this pressure corresponds to a Nonlinear Schr\"odinger Equation with nonlinear self-interaction potential obeying a power law (e.g the cubic defocusing NLS).\\
In the paper \cite{Na} regarding hydrodynamic nucleation of quantized vortex pairs in a polariton quantum fluid, theoretical predictions of different flow regimes are reported, to be depending on the Mach number, in particular the authors show that the low Mach number regime prevents the onset of the vortex nucleation mechanism.  

It is well known that a  way to obtain the incompressible system from the compressible one is to perform the so called {\em incompressible} or {low Mach number limit}. In fact, if we denote by $\e$ the Mach number 
$$\e=Mach\  number=\frac{typical \ fluid \ speed}{sound \ speed}$$
it makes sense to consider the limit $\e\to 0$.  When this situation occurs we observe that  the pressure becomes nearly constant and the fluid cannot generate density variations, so it behaves as an incompressible fluid. In order to study this dynamics on the  system \eqref{1}-\eqref{2} we perform the incompressible scaling given by
\begin{equation}
\label{inscal}
\re(x,t)=\rho\left( y\e^{-2}, s\e^{-1}\right), \ \Je=\e^{-1}J\left(y \e^{-2}, s \e^{-1}\right).
\end{equation}
With the scaling (\ref{inscal}) the system \eqref{1}-\eqref{2} becomes
\begin{equation}
\begin{split}
\partial_{t}\re+\dive \Je&=0\\
\partial_{t}\Je +\dive\left(\frac{\Je\otimes \Je}{\re}\right)+\e^{-2}\nabla p(\re)&=\dive(\re\nabla^{2}\log \re).
\label{ms}
\end{split}
\end{equation}
\subsection{Statement of the main result}
The goal of this paper will be the study of the limiting behavior of the system \eqref{ms} as $\e\to 0$.  Before giving a precise description of the limiting behavior of our system \eqref{ms} we need to define the framework where we are  going to set up our problem. In this paper we will always assume that $t\geq 0$ and $x\in \T^{n}$, where $\T^{n}$ is the $n$-dimensional torus.
\subsubsection{Weak solutions}
\label{sweak}
To simplify our notations from now on we will denote by $\Psi^{\e}$ the renormalized pressure, namely
\begin{equation}
\Psi^{\e}=\sqrt{\frac{(\re)^{\gamma}-1-\gamma(\re-1)}{\e^{2}\gamma(\gamma-1)}}
\label{rp}
\end{equation}
and by 
\begin{equation}
\se={\e}^{-1}{(\re-1)}
\label{df}
\end{equation}
the density fluctuation.
A natural framework to deal with the system \eqref{ms} is given by the space of finite initial energy. In fact  the energy associated to the system \eqref{ms} is given by
\begin{equation}
E(t)=\frac{1}{2}\int_{\T^{n}}\left(|\Le(t)|^{2}+(\Psi^{\e})^{2}+|\nabla\sqrt{\re}|^{2}\right)dx,
\label{energy}
\end{equation}
where $\Le=\Je/\sqrt{\re}$. So it seems now natural to introduce the following definition of {\em weak solution}.
\begin{definition}
\label{dweak}
We say that $(\re,\Je)$ is a weak solution for the system \eqref{ms} with initial data $(\re_{0}, \Je_{0})$ if there exists locally integrable functions $\sqrt{\re}, \Le$, such that $\sqrt{\re}\in L^{2}_{loc}((0,T);H^{1}_{loc}(\T^{n}))$, $\Le\in L^{2}_{loc}((0,T);L^{2}_{loc}(\T^{n}))$ and by defining $\re=(\sqrt{\re})^{2}$, $\Je=\sqrt{\re}\Le$ the following integral identity hold for any test function $\varphi\in C^{\infty}_{0}([0,T]\times \T^{n})$, $\varphi(\cdot, T)=0$
\begin{equation}
\int_{0}^{T}\!\!\int_{\T^{n}}\left(\re\partial_{t}\varphi+\Je\cdot\nabla\varphi\right) dxdt+\int_{\T^{n}}\re_{0}\varphi(0)dx=0
\end{equation}
and for any test function $\psi\in C^{\infty}_{0}([0,T]\times \T^{n};\R^{3})$, $\psi(\cdot, T)=0$
\begin{equation}
\begin{split}
\int_{0}^{T}\!\!\int_{\T^{n}}\Big(\Je\partial_{t}\psi+\Le\otimes\Le&:\nabla\psi+\frac{1}{\e^{2}}p(\re)\dive\psi\\
+4\nabla\sqrt{\re}\otimes \nabla\sqrt{\re}:\nabla\psi&-\re\nabla\dive\psi dxdt\Big)+\int_{\T^{n}}\Je_{0}\psi(0)dx=0.
\end{split}
\end{equation}
\end{definition}
The existence of irrotational weak solutions, including vacuum states, for finite energy  large data which are obtained by an $H^{1}$ wave function via a Madelung transformation has been proved by Antonelli and Marcati in \cite{AM1}, \cite{AM2}, by using dispersive analysis, local smoothing effects and polar factorization methods.  Moreover, by means of the convex integration methods, in \cite{DFM15} it has been proved that the system admits on the torus infinitely many global in time weak solutions for any sufficiently smooth initial data including the case of a vanishing initial density - the vacuum zones.
Existence of smooth solutions away from vacuum was proved before by Li and Marcati \cite {LiM} for small perturbations of quantum subsonic steady states.
Related results concerning the dynamics of quantum hydrodynamics systems can also be found in  \cite{GMSU95}, \cite{J09}, \cite{JHMS03}.
 
 \subsubsection{Initial data}
\label{Sinitial}
Since it is natural to work with weak solutions, that have bounded energy \eqref{energy} it is quite obvious to require that the initial data satisfy the following condition
\begin{equation}
E(0)=\frac{1}{2}\int_{\T^{n}}\left(|\Le_{0}|^{2}+(\Psi^{\e}_{0})^{2}+|\nabla\sqrt{\re_{0}}|^{2}\right)dx<+\infty,
\label{initial-energy}
\end{equation}
where 
$$\re|_{t=0}=\re_{0}\quad \Je|_{t=0}=\Je_{0}\quad \Le_{0}=\Je_{0}/\sqrt{\re_{0}}.$$

We perform our analysis by  considering sufficiently general initial data, that in a weak sense  can be called ill-prepared initial data since   we do not assume $\re_{0}=1$ and $\dive \Le_{0}=0$, but we simply require that
\begin{equation}
\Le_{0}\to \tilde{v}_{0} \qquad \mbox{strongly in} \,\, L^2(\T^n), 
\label{i1}
\end{equation}
\begin{equation}
\Psi^{\e}_{0}\to \delta_{0} \qquad\mbox{strongly in} \,\, L^2(\T^n), 
\label{i2}
\end{equation}
\begin{equation}
\sqrt{\re_{0}}-1\to 0 \qquad\mbox{strongly in} \,\, H^1(\T^n).
\label{i3}
\end{equation}
In particular \eqref{i2} implies that 
\begin{equation}
\se_{0}\to \tilde{\delta_{0}} \qquad\mbox{strongly in} \,\, L^{\gamma}(\T^n).
\label{i4}
\end{equation}
In order to apply in the sequel the relative entropy method we need some additional  regularity assumptions: $\sigma_{0}\in H^{s}(\T^{n})$, $\tilde{v}_{0}\in H^{s}(\T^{n})$, $P\tilde{v}_{0}\in H^{s}(\T^{n})$,
$Q\tilde{v}_{0}\in H^{s-1}(\T^{n})$, for any $s\geq n/2+1$, where $P$ and $Q$ denote the Leray projectors on the divergence free vector fields and the gradient vector fields respectively. 
\subsubsection{The limiting system}
As already said the aim of this paper is to perform the low Mach number limit for the system \eqref{ms}. If we  look at the second equation of \eqref{ms} (linear momentum equation) we can deduce that  as $\e\to 0$,  $\re$ behaves like $\tilde{\rho}+\e^{2}$ (where $\tilde{\rho}$ is a constant which by a simple scaling can always be assumed to be as $\tilde{\rho}=1$). So, at a formal level, we can see that as $\e\to 0$, the density $\re$ becomes constant, $\Le$ converges  to a soleinoidal vector field $v$. Hence   we end up with the following incompressible Euler system,
\begin{equation}
\begin{cases}
\dive v=0,\\
\partial_{t} v+\dive(v\otimes v)+\nabla\pi=0\\
v(x,0)=P\tilde{v}_{0}={v}_{0}.
\end{cases}
\label{euler}
\end{equation}
It is worthwhile, at this point, to recall the following classical result on the existence of  regular solutions for   the incompressible Euler  system \eqref{euler}, see Kato \cite{K72} and Lions \cite{LPL96}.

\begin{proposition}
\label{peuler}
Let the initial velocity field satisfy $ v_0 \in H^s (or \,\,  H^{s+1})$, $s \geq \frac{n}{2} +2$ with $\dive v_{0}=0.$
Then, there exists $0 < T^{\ast} < \infty$, the maximal existence time, and a unique smooth solution $(v, \pi)$ of the incompressible Euler equation\eqref{euler} on 
$[0, T^{\ast})$ with initial data $v_{0}$, satisfying for any $T < T^{\ast}$
$$\sup_{0\leq t< T^{\ast}}\left(\|v\|_{H^s} + \|\partial_t v\|
_{H^{s-1}} + \|\nabla\pi\|_{H^s} + \|\partial_t \nabla \pi\|_{H^{s-1}} \right)\leq M(T).$$
\end{proposition}

\subsubsection{Main result}
Now we are ready to state the main result we are going to prove in the paper
\begin{theorem}
\label{mt}
Assume that $(\re, \Le)$ is a weak solution of the quantum hydrodynamic system \eqref{ms}, in the sense of the Definition \ref{dweak}   and that initial data verify the conditions of Section \ref{Sinitial}. Let $T^{\ast}$ and $T^{\ast\ast}$ be as in the Propositions \ref{peuler}, \ref{poscill} respectively, then as $\e\to 0$ and for all $T<\min(T^{\ast}, T^{\ast\ast})$,  we have
\begin{itemize}
\item[i)] $\Le \rightharpoonup v$ weakly in $L^{\infty}(0,T;L^{2}(\T^{n}))$,
\item[ii)] $P\Le \rightarrow v$ strongly in $L^{\infty}(0,T;L^{2}(\T^{n}))$,
\end{itemize}
where $v$ is the unique local in time solution of the Euler system \eqref{euler} ($v\in L^{\infty}_{loc}(0, T^{\ast}, H^{s}(\T^{n}))$, $s > \frac{n}{2} +2$).
\end{theorem}

\subsection{Plan of the paper}

The low Mach number limit  for fluid dynamic models has been studied by many authors. See for examples the paper by P.L. Lions and Masmoudi \cite{L-P.L.M98}, Desjardin, Grenier, Lions, Masmoudi  \cite{DGLM99}, Desjardin Grenier \cite{DG99} for the case of the compressible Navier Stokes equations. The mathematical analysis is completely different in the case of Òwell-preparedÓ initial data ($\re_{0}=1$, $\dive(\Le_{0}/\sqrt{\re})=0$) or in the case of ``ill-prepared'' data. In the latter case the fluctuation of the fluid density are of the same order of the Mach number and the gradient part of the velocity develops fast time oscillations.
In fact the main issue in treating this kind of limits is the presence of acoustic waves  which propagate with high speed of order $1/\e$ and are supported by the gradient part of the velocity field. The main consequence is the loss of compactness of the  velocity field or of the momentum and the impossibility to define the limit of nonlinear quantities such as the convective term.
 The analysis at this point is different according to the space domain of the problem. We can say that in the case of the incompressible limit   the acoustic waves in general are well described  by a wave equation with a source term bounded in some suitable space.  Then, in the case of  an unbounded domain (whole or exterior domain) we  can observe that the acoustic waves redistribute their energy in the space and so one can exploit the dispersive properties of these waves to get the local  decay of the acoustic energy  and to recover compactness in time, see for example \cite{DG99}, \cite{DFN10}, \cite{DFN12}, \cite{DM12a}. In the case of a periodic domain we don't have a dispersion phenomenon but the waves interact with each other, so  in the spirit of Schochet \cite{Sc94} and \cite{Sc05} one has to introduce an operator which describes the oscillations in time so that they can be included in the energy estimates, see \cite{Ma01a}, \cite{Ma01b}, \cite{Donatelli-Trivisa-2008}, \cite{LC05}. 
 
 In this paper we will study  the incompressible limit in a periodic domain for the system of quantum hydrodynamics \eqref{ms} and, as explained above, the main issue is to control the time oscillations of the density fluctuation and of the momentum $\Je$.  In the Section 2, we start by recovering the standard energy estimates satisfied by the weak solutions of \eqref{ms}. Then, in Section 3 we introduce the  operator $\mathcal{L}$ which describes the time oscillations and will give a careful analysis of all its properties. In Section 4, in order to include the time oscillation in the energy and to study the convergence of our sequences, we introduce the relative   entropy functional $H^{\e}(t)$. This functional computes the error terms due to the fast time oscillations, namely the difference of our sequence and the limit solutions. In an heuristic way we can say that the entropy measures the error that we have when we pass from the  weak  to the strong convergence. We will be able to show that as $\e\to0$ the entropy goes to zero and so we get  the strong convergence of our sequences. This will lead to the proof of the main result in Section 5.
 
 For completeness we conclude this section by mentioning that in the same framework of the incompressible limits and its related problems and techniques  fits also the so called quasineutral limit in plasma physics or the zero electron mass limit, see for example  \cite{CDM13},  \cite{DFN15}, \cite{DM08}, \cite{DM12},  \cite{DM15}, \cite{JW06}, \cite{JLW08}, \cite{W04}, \cite{WJ06},

\section{Energy inequality and its consequences}
Taking into account the existence result of the Section \ref{sweak} we know that the weak solutions of the system \eqref{ms} satisfy the following energy bound
\begin{equation}
E(t)=\frac{1}{2}\int_{\T^{n}}\left(|\Le(t)|^{2}+(\Psi^{\e})^{2}+|\nabla\sqrt{\re}|^{2}\right)dx\leq E(0).
\label{bounds}
\end{equation}
Hence, by virtue of \eqref{initial-energy} and by the convexity of the function  $z \to z^{\gamma} -1 - \gamma (z -1)$ for $z \geq 0,$  the following convergence (up to a subsequence) hold:

\begin{equation}
 \re -1 \rightarrow 0 \quad \mbox{strongly in $L^{\infty}([0,T]; L^{\gamma}(\T^n)),$} 
 \label{c1}
\end{equation}

\begin{equation}
\Le\rightharpoonup v \qquad \mbox{weakly in $L^{\infty}([0,T]; L^2(\T^n))$.}  
\label{c2}
\end{equation}
Moreover since 
$$|\sqrt{z}-1|^{2}\leq M|z-1|^{\gamma}, \ |z-1|\geq \eta,\ \gamma\geq 1,$$
$$|\sqrt{z}-1|^{2}\leq M|z-1|^{2},\  z \geq 0,$$
from \eqref{c1} we have

\begin{equation}
 \sqrt{\re} -1 \rightarrow 0 \quad \mbox{strongly in $L^{\infty}([0,T]; L^{2}(\T^n)).$} 
 \label{c3}
\end{equation}
By rewriting the continuity equation $\eqref{ms}_{1}$ in the following way
$$\partial_{t}(\re-1)+\dive((\sqrt{\re}-1)\Le)+\dive\Le=0,$$
as $\e\to 0$ we infer that $v(x,t)$ is a divergence free vector field.
Unfortunately the previous convergences are not enough to pass into the limit in the system \eqref{ms} since we still do not  control the oscillations in time. This will be argument of the  next session.

\section{Study of the time oscillations}
In this section we try to understand the behavior of the oscillations in time in order to  prove that  they don't affect the limit system. 
In order to describe the time oscillations (following \cite{Sc05}) we introduce the group $\mathcal{L}(\tau)$, $\tau\in \R$, defined by $e^{\tau L}$, where $L$ is the operator on the space $\mathcal{D}'_{0}\times \mathcal{D}'$, where  $\mathcal{D}'_{0}=\{\phi\in \mathcal{D}' \mid \int \phi=0\}$ given by
\begin{equation}
L\begin{pmatrix}
     \phi \\
      v  
\end{pmatrix}=-
\begin{pmatrix}
\dive v\\
\nabla\phi
\end{pmatrix}. 
\label{o1}
\end{equation}
Notice that $\mathcal{L}$ is an isometry in  any $H^{s}$ space, $s\in \R$.
Now we introduce the following notations that we are going to use later on,
\begin{equation}
U^{\e}=\begin{pmatrix}
      \se\\
      Q(\Je)   \\
      \end{pmatrix},\qquad V^{\e}=\mathcal{L}\left(-\frac{t}{\e}\right)U^{\e},
\label{o2}
\end{equation}
 \begin{equation}
\bar{U}^{\e}=\begin{pmatrix}
      \Psi^{\e}\\
      Q(\Le)   \\
      \end{pmatrix},\qquad \bar{V}^{\e}=\mathcal{L}\left(-\frac{t}{\e}\right)\bar{U}^{\e},
\label{o3}
\end{equation}
and the following approximation holds
\begin{equation}
\|U^{\e}-\bar{U}^{\e}\|_{L^{\infty}(0,T;L^{2\gamma/(\gamma+1)}(\T^{n}))}\longrightarrow 0, \quad \mbox{as $\e\to 0$.}
\label{oo}
\end{equation}
 By using the notations \eqref{rp} and \eqref{df} we rewrite the system \eqref{ms} as follows
 \begin{equation}
\begin{split}
\e\partial_{t}\se+\dive Q(\Je)&=0\\
\e\partial_{t}Q(\Je) +\nabla\se&=\e G^{\e}
\label{mso}
\end{split}
\end{equation}
 where
 \begin{equation}
 G^{\e}=-Q\left[\dive\left(\frac{\Je\otimes \Je}{\re}\right)-\dive(\re\nabla^{2}\log \re)\right]-(\gamma-1)\nabla(\Psi^{\e})^{2}.
 \label{o4}
 \end{equation}
 By means of \eqref{o2}, the system \eqref{mso} has also the following form
 \begin{equation}
 \partial_{t} U^{\e}=\frac{1}{\e} LU^{\e}+\begin{pmatrix}0\\G^{\e}\end{pmatrix},
 \label{o5}
 \end{equation} 
 which is equivalent to
  \begin{equation}
 \partial_{t} V^{\e}=\mathcal{L}\left(-\frac{t}{\e}\right)\begin{pmatrix}0\\G^{\e}\end{pmatrix}.
 \label{o6}
 \end{equation} 
 From the energy bounds \eqref{bounds} we get  that $G^{\e}$ is bounded in $L^{2}([0,T];H^{-s}(\T^{n}))$,  $s>0$ uniformly in $\e$, hence $V^{\e}$ is compact in time  (the oscillations have been cancelled) and since $V^{\e}\in L^{\infty}([0,T];L^{\frac{2\gamma}{\gamma+1}}(\T^{n}))$, uniformly in $\e$, we get as $\e\to 0$,
 $$V^{\e}\rightarrow \bar{V}\ \text{strongly in $L^{p}([0,T];H^{-s'}(\T^{n}))$ for all $s'>s$ and $1<p<\infty$}.$$
 At this point it is important to remark that $\mathcal{L}\left(-\frac{t}{\e}\right)(0,G^{\e})$ can be considered as an almost periodic function in $\tau=t/\e$ and computing its mean values yields the definitions of the following bilinear forms (see \cite{Ma01b}, \cite{Sc94}).
 
 \begin{definition}
 For all divergence free vector field $v\in L^{2}(\T^{n})$ and all $V=(\psi, \nabla q)\in  L^{2}(\T^{n})$, we define the following linear and bilinear symmetric forms in $V$
 \begin{equation}
 B_{1}(v,V)=\lim_{\tau\to\infty}\frac{1}{\tau}\int_{0}^{\tau}\!\!\mathcal{L}(-s)\begin{pmatrix}0\\ \dive(v\otimes \mathcal{L}_{2}(s)V+\mathcal{L}_{2}(s)V\otimes v)\end{pmatrix}ds,
 \label{b1}
 \end{equation}
 and 
 \begin{equation}
 B_{2}(V,V)=\lim_{\tau\to\infty}\frac{1}{\tau}\int_{0}^{\tau}\!\!\mathcal{L}(-s)\begin{pmatrix}0\\ \dive(\mathcal{L}_{2}(s)V\otimes \mathcal{L}(s)V+(\gamma-1)\nabla(\mathcal{L}_{1}(s)V)^{2}\end{pmatrix}ds,
 \label{b2}
 \end{equation}
 \end{definition}
Now, in the same spirit as in \cite{Ma01b}, if we pass into the limit in \eqref{o6} we get that $\bar{V}$ satisfies the following equation
\begin{equation}
\partial_{t}\bar{V}+B_{1}(v, \bar{V})+B_{2}(\bar{V}, \bar{V})=0,
\label{o7}
\end{equation}
 where $B_{1}$ and $B_{2}$ are as in \eqref{b1} and \eqref{b2} respectively.
In a way similar to  \cite{Ma01b} we obtain the following local existence result for the system  \eqref{o7}.
 \begin{proposition}
 \label{poscill}
 Let us consider the following system 
 \begin{equation}
 \begin{cases}
 \partial_{t}V^{0}+B_{1}(v,V^{0})+B_{2}(V^{0}, V^{0})=0,\\
 V^{0}(0)=(\sigma_{0}, Q\tilde{v}_{0}),
 \end{cases}
 \label{o8}
 \end{equation}
 where $v$ is the solution of the incompressible Euler  problem \eqref{euler} and $(\sigma_{0}, Q\tilde{v}_{0})$ satisfy the regularity conditions of the Section \ref{Sinitial}. Then, there exists a maximal  existence time $0 < T^{\ast\ast} < \infty$ and a  unique local strong solution $V^{0}$  of   \eqref{o8} such that $V^{0}\in L^{\infty}([0,T^{\ast\ast});H^{s-1}(\T^{n}))\cap L^{2}([0,T^{\ast\ast});H^{s}(\T^{n}))$, for any  $s\geq n/2+1$.
 \end{proposition}
 
 \begin{remark}
 It is important to notice that in the case of well prepared initial data, i.e. $V^{0}(0)=(\sigma_{0}, Q\tilde{v}_{0})=0$, the solution of the system \eqref{o8} is given by $V^{0}=0$. This means that the oscillations with respect to   time  vanish  and so $\Le\rightarrow v$ strongly in $L^{\infty}([0,T], L^{2}(\T^{n}))$. But, for the general initial data, since  the oscillation part with
respect to time $t/\e$ does not vanish, there are oscillations in time of the solution sequence. 
  \end{remark}
Now we report three technical proposition concerning the properties of the linear and bilinear forms $B_{1}$ and $B_{2}$ that we will use in the sequel. For their proofs  we refer to \cite{Ma01b}.
\begin{proposition}
For all $v, V, V_{1}, V_{2}$ we have
\begin{equation}
\int B_{1}(v,V)V=0 \qquad \mbox{and} \qquad \int B_{1}(v,V_{1})V_{2}+ B_{1}(v,V_{2})V_{1}=0 
\label{o9}
\end{equation}
\begin{equation}
 \int B_{2}(V,V)V=0
\quad \mbox{and} \quad \int B_{2}(V_{1},V_{1})V_{2}+2 B_{2}(V_{1},V_{2})V_{1}=0.
\label{o10}
\end{equation}
\end{proposition}

 \begin{proposition}
 For all $v\in L^{p}(0,T;L^{2}(\T^{n}))$ and $V\in L^{q}(0,T;L^{2}(\T^{n}))$, as $\e\to 0$ we have the following weak convergence ($p$ and $q$ are such that the products are well defined) 

 \begin{equation}
 \mathcal{L}\left(-\frac{t}{\e}\right)\begin{pmatrix}0\\ \dive(v\otimes \mathcal{L}_{2}\left(\frac{t}{\e}\right)V+\mathcal{L}_{2}\left(\frac{t}{\e}\right)V\otimes v)\end{pmatrix}\xrightarrow[weakly]\quad B_{1}(v,V) 
 \label{pp1}
 \end{equation}

\begin{equation}
 \mathcal{L}\left(-\frac{t}{\e}\right)\begin{pmatrix}0\\ \dive(\mathcal{L}_{2}\left(\frac{t}{\e}\right)V\otimes \mathcal{L}_{2}\left(\frac{t}{\e}\right)V)+(\gamma-1)\nabla(\mathcal{L}_{1}\left(\frac{t}{\e}\right)V)^{2}\end{pmatrix} \xrightarrow[weakly]\quad B_{2}(V,V).
 \label{pp2}
 \end{equation}
 \end{proposition}

\begin{proposition}
For any $V_{1}\in L^{q}(0,T;H^{s}(\T^{n}))$ and $V_{2}\in L^{p}(0,T;H^{-s}(\T^{n}))$, with $s\in \R$, $1/p+1/q=1$ one has as $\e\to 0$
 \begin{equation}
 \begin{split}
\mathcal{L}\left(-\frac{t}{\e}\right)&\begin{pmatrix}0\\ \dive(\mathcal{L}_{2}\left(\frac{t}{\e}\right)V_{1}\otimes \mathcal{L}_{2}\left(\frac{t}{\e}\right)V_{2}+
 \mathcal{L}_{2}\left(\frac{t}{\e}\right)V_{2}\otimes \mathcal{L}_{2}\left(\frac{t}{\e}\right)V_{1}\end{pmatrix}\\
 &+\mathcal{L}\left(-\frac{t}{\e}\right)\begin{pmatrix}0\\ (\gamma-1)\nabla(\mathcal{L}_{1}\left(\frac{t}{\e}\right)V_{1}\mathcal{L}_{1}\left(\frac{t}{\e}\right)V_{2})\end{pmatrix} \xrightarrow[weakly]\quad B_{2}(V_{1},V_{2}).
 \end{split}
 \label{pp3}
 \end{equation}
 It is also possible to extend \eqref{pp3} to the case where we replace $V_{2}$ in the left hand side by a sequence $V_{2}^{\e}$ that converges strongly to $V_{2}$ in $L^{p}(0,T;H^{-s}(\T^{n}))$.

\end{proposition}

\section{Relative entropy}

In order to prove the convergence stated in the Theorem \ref{mt} we introduce the following relative entropy functional
\begin{equation}
\mathcal{H}^{\e}(t)=\frac{1}{2}\int_{\T^{n}}\left\{|\Le -v- \mathcal{L}_{2}\left(\frac{t}{\e}\right)V^{0}|^{2}+|\Psi^{\e}-\mathcal{L}_{1}\left(\frac{t}{\e}\right)V^{0}|^{2}+|\nabla\sqrt{\re}|^{2}\right\}dx.
\label{entropy}
\end{equation}
The entropy describes the difference between the solutions of the scaled quantum hydrodynamic system \eqref{ms} and the limit solution, namely the solution $v$ of the Euler system \eqref{euler} and the fast time oscillations. 
The goal of this section will be to recover uniform estimates in $\e$ for $\mathcal{H}^{\e}(t)$ and to show that the relative entropy vanishes as $\e\to 0$, yielding the strong convergence of our solutions. First of all we recall that the solutions of \eqref{ms} satisfy the following energy bound
\begin{equation}
\begin{split}
\frac{1}{2}\int_{\T^{n}}\big(|\Le(t)|^{2}+(\Psi^{\e})^{2}&+|\nabla\sqrt{\re}|^{2}\big)dx\\&\leq\frac{1}{2}\int_{\T^{n}}\left(|\Le_{0}|^{2}+(\Psi^{\e}_{0})^{2}+|\nabla\sqrt{\re_{0}}|^{2}\right)dx.
\end{split}
\label{h1}
\end{equation}
Moreover the solution $v$  of the Euler system \eqref{euler}  satisfies the conservation of energy
\begin{equation}
\frac{1}{2}\int_{\T^{n}}|v|^{2}dx=\frac{1}{2}\int_{\T^{n}}|v_{0}|^{2}dx,
\label{h2}
\end{equation}
while the solution  $V^{0}$ of \eqref{o8}, taking into account \eqref{o9} and \eqref{o10}, satisfies
\begin{equation}
\frac{1}{2}\int_{\T^{n}}|V^{0}|^{2}dx=\frac{1}{2}\int_{\T^{n}}(|\sigma_{0}|^{2}+|Q\tilde{v}_{0}|^{2})dx,
\label{h3}
\end{equation}
Now if we use $\mathcal{L}_{1}\left(\frac{t}{\e}\right)V^{0}$ as test function for the weak formulation of the mass conservation equation $\eqref{ms}_{1}$, we have
\begin{align}
\int_{\T^{n}}\!\!&\mathcal{L}_{1}\left(\frac{t}{\e}\right)V^{0}\se dx-\int_{0}^{t}\!\! \int_{\T^{n}}\mathcal{L}_{1}\left(\frac{s}{\e}\right)\partial_{s}V^{0}\se dxds \nonumber\\
&+\frac{1}{\e}\int_{0}^{t}\!\! \int_{\T^{n}}\!\!\left(\!\dive(\sqrt{\re}\Le)\mathcal{L}_{1}\left(\frac{s}{\e}\right)V^{0}+\dive\!\!\left(\mathcal{L}_{2}\left(\frac{s}{\e}\right)V^{0}\right)\se \!\!\right)dxds\nonumber\\
&=\int_{\T^{n}}\sigma_{0}\se(0) dx.
\label{h4}
\end{align}
By using $v$ and then  $\mathcal{L}_{2}\left(\frac{t}{\e}\right)V^{0}$ as test functions in the momentum equation $\eqref{ms}_{2}$, we have
\begin{align}
\int_{\T^{n}}\sqrt{\re}\Le vdx&+\int_{0}^{t}\!\! \int_{\T^{n}}\sqrt{\re}\Le(v\cdot\nabla v+\nabla\pi)dxds-\int_{0}^{t}\!\! \int_{\T^{n}}\Le\otimes\Le\cdot\nabla v dxds\nonumber\\
&+\int_{0}^{t}\!\! \int_{\T^{n}}\nabla \sqrt{\re}\otimes \nabla \sqrt{\re}:\nabla v dxds=\int_{\T^{n}}\sqrt{\re_{0}}\Le_{0}v_{0} dx,
\label{h5}
\end{align}
\begin{align}
&\int_{\T^{n}}\sqrt{\re}\Le \mathcal{L}_{2}\left(\frac{t}{\e}\right)V^{0}dx-\int_{0}^{t}\!\! \int_{\T^{n}}\Le\otimes\Le\cdot\nabla \mathcal{L}_{2}\left(\frac{s}{\e}\right)V^{0}dxds\nonumber\\
&+\int_{0}^{t}\!\! \int_{\T^{n}}4\nabla \sqrt{\re}\otimes \nabla \sqrt{\re}:\nabla \mathcal{L}_{2}\left(\frac{s}{\e}\right)V^{0}-\re\Delta\dive \mathcal{L}_{2}\left(\frac{s}{\e}\right)V^{0} dxds\nonumber\\
&-\int_{0}^{t}\!\! \int_{\T^{n}}\mathcal{L}_{2}\left(\frac{s}{\e}\right)\partial_{s}V^{0}\sqrt{\re}\Le-\frac{1}{\e}\sqrt{\re}\Le\nabla \mathcal{L}_{1}\left(\frac{s}{\e}\right)V^{0}dxds\nonumber\\
&-\int_{0}^{t}\!\! \int_{\T^{n}}\left(\frac{1}{\e}\se+(\gamma-1)(\Psi^{\e})^{2}\right)\dive\left(\mathcal{L}_{2}\left(\frac{s}{\e}\right)V^{0}\right)dxds\nonumber\\
&=\int_{\T^{n}}\sqrt{\re_{0}}\Le_{0}Q\tilde{v}_{0} dx.
\label{h6}
\end{align}
Now taking into account that 
$$\int_{0}^{t}\!\! \int_{\T^{n}}\mathcal{L}\left(\frac{s}{\e}\right)\partial_{s}V^{0}U^{\e}dxds=\int_{0}^{t}\!\! \int_{\T^{n}}\partial_{s}V^{0}V^{\e}dxds,$$
we sum up \eqref{h1}, \eqref{h2}, \eqref{h3} and subtract \eqref{h4}, \eqref{h5}, \eqref{h6}, therefore we get that following inequality for $H^{\e}(t)$
\begin{equation}
\mathcal{H}^{\e}(t)\leq I^{\e}+A^{\e}+B^{\e}+C^{\e}.
\label{h7}
\end{equation}
Where we set
\begin{equation}
I^{\e}=\mathcal{H}^{\e}(0)+ \int_{\T^{n}}v\mathcal{L}_{2}\left(\frac{t}{\e}\right)V^{0}dx+ \int_{\T^{n}}(\sqrt{\re_{0}}-1)\Le_{0}(v_{0}+Q\tilde{v}_{0})dx
\label{h8}
\end{equation}
\begin{equation}
A^{\e}= -\int_{\T^{n}}(\Psi^{\e}-\se) \mathcal{L}_{1}\left(\frac{t}{\e}\right)V^{0}dx+ \int_{\T^{n}}(\sqrt{\re}-1)\Le(v+\mathcal{L}_{2}\left(\frac{t}{\e}\right)V^{0})dx,
\label{h9}
\end{equation}

\begin{align}
B^{\e}= &\int_{0}^{t}\!\! \int_{\T^{n}}\!\!\left(\sqrt{\re}\Le(v\cdot\nabla v+\nabla\pi)-\Le\otimes\Le\cdot\!\nabla\! \left(v+\mathcal{L}_{2}\left(\frac{s}{\e}\right)V^{0}\!\right)\!\!\right) dxds\nonumber\\
&+\int_{0}^{t}\!\! \int_{\T^{n}} (\gamma-1)(\Psi^{\e})^{2}\dive\left(\mathcal{L}_{2}\left(\frac{s}{\e}\right)V^{0}\right)dxds\nonumber\\
&+\int_{0}^{t}\!\! \int_{\T^{n}}(B_{1}(v, V^{0})V^{\e}+B_{1}(V^{0}, V^{0})V^{\e})dxds,
\label{h10}
\end{align}

\begin{align}
C^{\e}= &\int_{0}^{t}\!\! \int_{\T^{n}}4\nabla \sqrt{\re}\otimes \nabla \sqrt{\re}:(\nabla v +\nabla \mathcal{L}_{2}\left(\frac{s}{\e}\right)V^{0})dxds\nonumber\\
&-\int_{0}^{t}\!\! \int_{\T^{n}}\re\Delta\dive \mathcal{L}_{2}\left(\frac{s}{\e}\right)V^{0} dxds.
\label{h11}
\end{align}

\subsection{Uniform estimates for $\mathcal{H}^{\e}(t)$}
Here we will estimate uniformly in $\e$ the righthand side of \eqref{h7}. In what follows we will denote by $r^{\e}(t)$ any term such that $r^{\e}(t)\to 0$, as $\e\to 0$ and by $M(T)$ a constant that depends only on $T=\min(T^{\ast}, T^{\ast\ast})$. We start with $I^{\e}$. By taking into account the assumptions on the initial data of the Section \ref{Sinitial}, the properties $v$ solution of the system \eqref{euler} and of the operator $\mathcal{L}$ we get
\begin{equation}
| I^{\e}|\leq C\mathcal{H}^{\e}(0)+\e E(0)\|v_{0}+Q\tilde{v}_{0}\|_{H^{s}(\T^{n})}\leq \mathcal{H}^{\e}(0)+r^{\e}(t),
\label{h12}
\end{equation}
where $C>0$ is a constant.
In order to estimate $A^{\e}$ we use \eqref{oo} and the regularity of $V^{0}$ and we get
\begin{align}
|A^{\e}|&\leq r^{\e}(t)M(T)+\left|\int_{\T^{n}}(\Psi^{\e}-\se) \mathcal{L}_{1}\left(\frac{t}{\e}\right)V^{0}dx\right|\nonumber\\
&\leq r^{\e}(t)M(T)+M(T)\|U^{\e}-\bar{U}^{\e}\|_{L^{\infty}(0,T;L^{2\gamma/\gamma+1})}\leq  r^{\e}(t)M(T) .
\label{h13}
\end{align}
In the same spirit we estimate $C^{\e}$ and we end up with
\begin{align}
|C^{\e}|&\leq M(T)\int_{0}^{t}\|\nabla \sqrt{\re}\|_{L^{2}(\T^{n})}ds+2\int_{0}^{t}\!\! \int_{\T^{n}}( \sqrt{\re}-1)\nabla \sqrt{\re}\nabla \mathcal{L}_{2}\left(\frac{s}{\e}\right)V^{0}dxds\nonumber\\
&\ +2\int_{0}^{t}\!\! \int_{\T^{n}}\nabla  \sqrt{\re}\nabla \mathcal{L}_{2}\left(\frac{s}{\e}\right)V^{0}dxds\nonumber\\
&\leq M(T)\int_{0}^{t}\mathcal{H}^{\e}(s)ds.
\label{h14}
\end{align}
The term $B^{\e}$ deserves some more attention, first of all we split it in three parts as follows,
\begin{equation}
|B^{\e}|\leq |B^{\e}_{1}|+|B^{\e}_{2}|+|B^{\e}_{3}|,
\label{h15}
\end{equation}
and then we estimate  each one of the  three parts. We start with $B^{\e}_{1}$,
\begin{equation}
|B^{\e}_{1}|=\left|\int_{0}^{t}\!\! \int_{\T^{n}}\sqrt{\re}\Le\nabla\pi dxds\right|\leq  r^{\e}(t)+M(T)\int_{0}^{t}\mathcal{H}^{\e}(s)ds.
\label{h16}
\end{equation}
Then, as $\e\to 0$ we also have,

\begin{align}
|B^{\e}_{2}|&=\left|\int_{0}^{t}\!\! \int_{\T^{n}}(B_{1}(v, V^{0})V^{\e}+B_{2}(V^{0}, V^{0})V^{\e})dxds\right|\nonumber\\
&\leq\int_{0}^{t}\!\! \int_{\T^{n}}\left|(B_{1}(v, V^{0})\bar{V}+B_{2}(V^{0}, V^{0})\bar{V})\right|dxds+ r^{\e}(t),
\label{h17}
\end{align}

\begin{align}
|B^{\e}_{3}|&=\Bigg|\int_{0}^{t}\!\! \int_{\T^{n}}\sqrt{\re}\Le (v\cdot\nabla) vdxds-\int_{0}^{t}\!\! \int_{\T^{n}}\Le\otimes\Le\nabla \left(v+\mathcal{L}_{2}\left(\frac{s}{\e}\right)V^{0}\right)dxds\nonumber\\
&\quad-\int_{0}^{t}\!\! \int_{\T^{n}} (\gamma-1)(\Psi^{\e})^{2}\dive\left(\mathcal{L}_{2}\left(\frac{s}{\e}\right)V^{0}\right)dxds\Bigg|\nonumber\\
&\leq\Bigg|-\int_{0}^{t}\!\! \int_{\T^{n}}(\Le -v- \mathcal{L}_{2}\left(\frac{s}{\e}\right)V^{0})\otimes (\Le -v- \mathcal{L}_{2}\left(\frac{s}{\e}\right)V^{0})\nonumber\\
&\qquad\cdot\nabla(v+ \mathcal{L}\left(\frac{s}{\e}\right)V^{0})dxds\nonumber\\
&\quad-(\gamma-1)\int_{0}^{t}\!\! \int_{\T^{n}}|\Psi^{\e}-\mathcal{L}_{1}\left(\frac{s}{\e}\right)V^{0}|^{2}\dive\left(\mathcal{L}_{2}\left(\frac{s}{\e}\right)V^{0}\right)dxds\nonumber\\
&\quad-\int_{0}^{t}\!\! \int_{\T^{n}}\Big [\Le\otimes\left(v+\mathcal{L}_{2}\left(\frac{s}{\e}\right)V^{0}\right)+
\left(v+\mathcal{L}_{2}\left(\frac{s}{\e}\right)V^{0}\right)\otimes\Le\Big ]\nonumber\\
&\qquad\cdot\nabla\left(v+ \mathcal{L}_{2}\left(\frac{s}{\e}\right)V^{0}\right)dxds\nonumber\\
&\quad+\int_{0}^{t}\!\! \int_{\T^{n}}\left(v+\mathcal{L}_{2}\left(\frac{s}{\e}\right)V^{0}\right)\otimes\left(v+\mathcal{L}_{2}\left(\frac{s}{\e}\right)V^{0}\right)\nonumber\\
&\qquad\cdot\nabla\left(v+\mathcal{L}_{2}\left(\frac{s}{\e}\right)V^{0}\right)dxds\nonumber\\
&\quad+\int_{0}^{t}\!\! \int_{\T^{n}}\Big\{(\gamma-1)\left|\mathcal{L}_{1}\left(\frac{s}{\e}\right)V^{0}\right|^{2}\dive\left(\mathcal{L}_{2}\left(\frac{s}{\e}\right)V^{0}\right)\nonumber\\
&\quad-(\gamma-1)\mathcal{L}_{1}\left(\frac{s}{\e}\right)V^{0}\Psi^{\e}\dive\left(\mathcal{L}_{2}\left(\frac{t}{\e}\right)V^{0}\right)\Big\}dxds\Bigg|+r^{\e}(t).
\label{h18}
\end{align}
Now by using \eqref{pp2} we have
\begin{align}
\int_{0}^{t}\!\! &\int_{\T^{n}}\left(\mathcal{L}_{2}\left(\frac{s}{\e}\right)V^{0}\right)\otimes\left(\mathcal{L}_{2}\left(\frac{s}{\e}\right)V^{0}\right)\cdot\nabla\left(v+\mathcal{L}_{2}\left(\frac{s}{\e}\right)V^{0}\right)dxds\nonumber\\
&(\gamma-1)\int_{0}^{t}\!\! \int_{\T^{n}}\left|\mathcal{L}_{1}\left(\frac{s}{\e}\right)V^{0}\right|^{2}\dive\left(\mathcal{L}_{2}\left(\frac{s}{\e}\right)V^{0}\right)dxds\nonumber\\
&=-\int_{0}^{t}\!\!\int_{\T^{n}}\Big[\dive\left(\mathcal{L}_{2}\left(\frac{s}{\e}\right)V^{0}\right)\otimes\left(\mathcal{L}_{2}\left(\frac{s}{\e}\right)V^{0}\right)\nonumber\\
&+(\gamma-1)\nabla\left|\mathcal{L}_{1}\left(\frac{s}{\e}\right)V^{0}\right|^{2}\Big]\cdot\left(V^{0}+\begin{pmatrix}0\\v\end{pmatrix}\right)dxds\nonumber\\
&=-\int_{0}^{t}\!\!\int_{\T^{n}}\mathcal{L}\left(\frac{s}{\e}\right)V^{0}\begin{pmatrix}0\\ 
\dive\left(\mathcal{L}_{2}\left(\frac{s}{\e}\right)V^{0}\right)\otimes\left(\mathcal{L}_{2}\left(\frac{s}{\e}\right)V^{0}\right)+(\gamma-1)\nabla\left|\mathcal{L}_{1}\left(\frac{s}{\e}\right)V^{0}\right|^{2}
\end{pmatrix}\nonumber\\
&\cdot\left(V^{0}+\begin{pmatrix}0\\v\end{pmatrix}\right)dxds-\int_{0}^{t}\!\!\int_{\T^{n}}B_{2}(V^{0},V^{0})\cdot\left(V^{0}+\begin{pmatrix}0\\v\end{pmatrix}\right)dxds+r^{\e}(t)=r^{\e}(t).
\label{h19}
\end{align}
In the same way if we use \eqref{pp3} we get
\begin{align}
&-\int_{0}^{t}\!\!\int_{\T^{n}}\Big[\mathcal{L}_{2}\left(\frac{s}{\e}\right)V^{0}\otimes\Le+\Le\otimes \mathcal{L}_{2}\left(\frac{s}{\e}\right)V^{0}\Big]\cdot \nabla \left(v+\mathcal{L}_{2}\left(\frac{s}{\e}\right)V^{0}\right)dxds\nonumber\\
&-\int_{0}^{t}\!\!\int_{\T^{n}}(\gamma-1)\mathcal{L}_{1}\left(\frac{s}{\e}\right)V^{0}\Psi^{\e}\dive\left(\mathcal{L}_{2}\left(\frac{s}{\e}\right)V^{0}\right)dxds\nonumber\\
&=\int_{0}^{t}\!\!\int_{\T^{n}}\Big[\dive\left(\mathcal{L}_{2}\left(\frac{s}{\e}\right)V^{0}\otimes\Le+\Le\otimes \mathcal{L}_{2}\left(\frac{s}{\e}\right)V^{0}\right)\nonumber\\
&+(\gamma-1)\nabla\left(\mathcal{L}_{1}\left(\frac{s}{\e}\right)V^{0}\Psi^{\e}\right)\cdot \left(v+\mathcal{L}_{2}\left(\frac{s}{\e}\right)V^{0}\right)dxds\nonumber\\
&=\int_{0}^{t}\!\!\int_{\T^{n}}\mathcal{L}\left(\frac{s}{\e}\right)\begin{pmatrix}0\\
\dive\left(\mathcal{L}_{2}\left(\frac{s}{\e}\right)V^{0}\otimes\Le+\Le\otimes \mathcal{L}_{2}\left(\frac{s}{\e}\right)V^{0}\right)+(\gamma-1)\nabla\left(\mathcal{L}_{1}\left(\frac{s}{\e}\right)V^{0}\Psi^{\e}\right)
\end{pmatrix}\nonumber\\
&\cdot\left(V^{0}+\begin{pmatrix}0\\v\end{pmatrix}\right)dxds\nonumber\\
&=\int_{0}^{t}\!\!\int_{\T^{n}}\left(2B_{2}(V^{0},\bar{V})+B_{1}(v,V^{0})\right)\cdot\left(V^{0}+\begin{pmatrix}0\\v\end{pmatrix}\right)dxds+r^{\e}(t)\nonumber\\
&=\int_{0}^{t}\!\!\int_{\T^{n}}2B_{2}(V^{0},\bar{V})(v,V^{0})\cdot\left(V^{0}+\begin{pmatrix}0\\v\end{pmatrix}\right)dxds+r^{\e}(t).
\label{h20}
\end{align}
By standard computations we also get
\begin{align}
\int_{0}^{t}\!\!\int_{\T^{n}}&\dive(v\otimes \Le+\Le\otimes v)\cdot\left(v+ \mathcal{L}\left(\frac{s}{\e}\right)V^{0}\right)dxds\nonumber\\
&=\int_{0}^{t}\!\!\int_{\T^{n}}B_{1}(v,\bar{V})V^{0}dxds+r^{\e}(t).
\label{h21}
\end{align}
 and
\begin{align}
\int_{0}^{t}\!\!\int_{\T^{n}}&\dive(v\otimes \mathcal{L}_{2}\left(\frac{s}{\e}\right)V^{0}+\mathcal{L}_{2}\left(\frac{s}{\e}\right)V^{0}\otimes v)\cdot\left(v+ \mathcal{L}\left(\frac{s}{\e}\right)V^{0}\right)dxds\nonumber\\
&=\int_{0}^{t}\!\!\int_{\T^{n}}B_{1}(v,V^{0})V^{0}dxds+r^{\e}(t)=r^{\e}(t).
\label{h22}
\end{align} 
By adding up \eqref{h16}-\eqref{h15} and by using the properties \eqref{o9} and \eqref{o10}, the term \eqref{h15} assumes the form
\begin{align}
|B^{\e}|&\leq r^{\e}(t)+M(T)\int_{0}^{t}H^{\e}(\tau)d\tau\nonumber\\
&\quad+\int_{0}^{t}\!\!\int_{\T^{n}}(2B_{2}(V^{0},\bar{V})V^{0}+B_{2}(V^{0},V^{0})\bar{V})dxds\nonumber\\
&\quad+\!\!\int_{\T^{n}}(B_{1}(v,\bar{V})V^{0}+B_{1}(v,V^{0})\bar{V})dxds\nonumber\\
&=r^{\e}(t)+M(T)\int_{0}^{t}H^{\e}(s)ds.
\label{h23}
\end{align}
By considering \eqref{h13}, \eqref{h14} and \eqref{h23} together, we can conclude that the relative entropy $\mathcal{H}^{\e}(t)$ satisfies the following inequality,
\begin{equation}
\mathcal{H}^{\e}(t)\leq C\mathcal{H}^{\e}(0) + M(T)\int_{0}^{t}\mathcal{H}^{\e}(s)ds+r^{\e}(t),
\label{h24}
\end{equation}
from which,  since $r^{\e}(t)\to 0$ as $\e\to 0$ and, by using Gronwall's inequality, we get there exists a constant $M>0$ such that,
\begin{equation}
\mathcal{H}^{\e}(t)\leq M,\qquad \mbox{for any $t\in [0,T]$, uniformely in $\e$.}
\label{h25}
\end{equation}
\subsection{Convergence of the relative entropy}
Because of the bound \eqref{h25} it makes sense to define the following quantity
$$\eta(t)=\limsup_{\e\to 0} \mathcal{H}^{\e}(t).$$
We get from \eqref{h24} that
\begin{equation}
\eta(t)\leq \eta(0) + M(T)\int_{0}^{t}\eta(s)ds.
\label{h26}
\end{equation}
Since the initial conditions \eqref{i1}-\eqref{i3} entail that $\eta(0)\equiv 0$, from \eqref{h26} we can conclude that 
\begin{equation}
\eta(t)=\limsup_{\e\to 0} \mathcal{H}^{\e}(t)=0 \qquad \mbox{for any $t\in [0,T]$}.
\label{h27}
\end{equation}

\section{Proof of the Theorem \ref{mt}}
Because of the previous estimates we have now the uniform "a priori" bounds needed to prove the Theorem \ref{mt}. Therefore we have that (i) is a consequence of \eqref{c2}, while  (ii) follows from \eqref{h27} and the following estimate
\begin{align*}
\sup_{0\leq t\leq T}\|P\Le-v\|_{L^{2}(\T^{n})}&=\sup_{0\leq t\leq T}\|P\left(\Le-v-\mathcal{L}_{2}\left(\frac{t}{\e}\right)V^{0}\right)\|_{L^{2}(\T^{n})}\\
&\leq \sup_{0\leq t\leq T}\|\Le-v-\mathcal{L}_{2}\left(\frac{t}{\e}\right)V^{0}\|_{L^{2}(\T^{n})}\\
&\leq \sup_{0\leq t\leq T} \mathcal{H}^{\e}(t)\rightarrow 0 \quad \mbox{as $\e\to 0$}.
\end{align*}

\bibliographystyle{amsplain}

\end{document}